\input amstex
\documentstyle{amsppt}
\document
\magnification=1200
\NoBlackBoxes
\nologo
\vsize16cm


\bigskip

\centerline{\bf ITERATED SHIMURA INTEGRALS}

\medskip

\centerline{\bf Yuri I. Manin}

\medskip

\centerline{\it Max--Planck--Institut f\"ur Mathematik, Bonn, Germany,}

\centerline{\it and Northwestern University, Evanston, USA}

\bigskip

{\bf Abstract.} In this paper I continue the study
of iterated integrals of modular forms  and noncommutative modular
symbols for $\Gamma \subset SL(2,\bold{Z})$ started in [Ma3].
Main new results involve a description of the iterated
Shimura cohomology and the image of the iterated Shimura cocycle class
inside it. The concluding section of the paper
contains a concise review of the classical modular
symbols for $SL(2)$ and a discussion of open problems.

\bigskip

\centerline{\bf \S 0. Introduction}

\medskip

Let $M$ be a linearly connected space, and $G$ a group
acting on it. Then $G$ acts on the fundamental
groupoid of $M$ thus creating a situation
where the well known formalism of cohomology of
$G$ with noncommutative coefficients applies.

\smallskip

If $M$ is a differentiable manifold, Chen iterated
integrals produce a representation of the fundamental groupoid
so that we get relations between such integrals reflecting
the action of $G$.

\smallskip

In [Ma3] I have studied this situation for the case
when $M$ is the upper complex half plane partially completed
by cusps, and the iterated integrals involve
cusp forms (and eventually Eisenstein series). The questions asked 
and the form of answers I would like to get in this case
were motivated by Drinfeld's associators and the classical theory of ordinary
integrals including the basics of Mellin transform and modular symbols.

\smallskip

Here I continue this study, stressing the Shimura approach
to the $SL(2)$--modular symbols of arbitrary weight
and attempting its iterated extension.

\smallskip

The paper is structured as follows. In \S 1
the notation and some background of noncommutative
group cohomology is reviewed. In \S 2 the theory
of the iterated Shimura cocycle is given. Finally,
\S 3 sketches the classical theory of modular symbols
and discusses open problems.
The reader might prefer to read this section first, as
a motivation for our attempt to produce its iterated version.

\newpage

\centerline{\bf \S 1. Noncommutative cohomology}

\smallskip

\centerline{\bf and abstract Shimura--Eichler relations}

\medskip

In this section I set notation and collect some
general background facts.

\medskip

{\bf 1.1. Noncommutative group cohomology: a general formalism.} Let $G$ be a group, and $N$ a group endowed with a left action of $G$ by automorphisms: $(g,n)\mapsto gn$.
Generally, both $G$ and $N$ can be noncommutative, and the group
laws are written multiplicatively. 

\smallskip

The set of 1--cocycles is defined by
$$
Z^1(G,N) := \{\, u:\,G\to N\,|\,u(g_1g_2)= u(g_1)\,g_1u(g_2)\,\} .
$$
It follows that $u(1_G)=1_N$.

\smallskip
Two cocycles are cohomological, $u^{\prime} \sim u$, iff
there exists an $n\in N$ such that for all $g\in G$ we have $u^{\prime}(g)=n^{-1}\,u(g)\cdot gn.$ This is an equivalence relation,
and by definition, 
$$
H^1(G,N) := Z^1(G,N)/(\sim ).
$$
This is a set with a marked point: the class of the
trivial cocycles  $u_n(g) =n^{-1} gn$.

\smallskip

Assume now that $G$ is embedded into a larger group $G\subset H$,
$[H:G] < \infty$. Denote by $N_H$ the induced noncommutative $H$--module:
$N_H$ is the space of $G$--covariant maps $\phi:\,H\to N,\,\phi (gh)=
g\phi (h)$ for all $g\in G$, $h\in H$, with
pointwise multiplication and the left action of $H$
$$
(h\phi)(h^{\prime}):=\phi (h^{\prime}h)\, .
$$
The map $N_H\to N$:\, $\phi \mapsto \phi (1_G)$, is a group homomorphism,
compatible with the action of $G$. Hence it induces
a map of pointed sets $Z^1(H,N_H)\to Z^1(G,N).$ One easily checks that 
cohomological cocycles go to the cohomological ones so that
we have an induced map $c:\,H^1(H,N_H)\to H^1(G,N).$ 

\smallskip

{\bf 1.1.1. Proposition (noncommutative Shapiro Lemma).} {\it
The map $c$ is a bijection.}

\smallskip

For a proof, see [PlRap], I.1.3. Here we only describe for future use
a map $Z^1(G,N) \to Z^1(H,N_H)$ which sends equivalent cocycles
to equivalent ones and induces the inverse map
$c^{-1}:\,H^1(G,N) \to H^1(H,N_H).$

\smallskip

To this end we will slghtly modify notation:
for a cocycle $u:\, G \to N$ and $g\in G$,
we will now denote by $u_g\in N$ the former $u(g)$. We want
to produce from $u$ a cocycle $\tilde{u}$ whose value
at $h\in H$ will be denoted $\tilde{u}_h \in N_H$. 
Thus, $\tilde{u}_h$ is a $G$--covariant function
$H\to N$ whose value at $h^{\prime}\in H$ will be denoted
$\tilde{u}_h (h^{\prime}).$ A well defined
prescription for obtaining this value, according to [PlRap],
requires a choice of representatives of $G\setminus H$ in $H$:
let $H=\coprod _i Gh_i$, such that $G$ is represented by $1_H$.
Then we put, for any $g, g^{\prime} \in G$:
$$
\tilde{u}_{gh_i}(g^{\prime}h_j):= 
g^{\prime}u_{g_{ji}}
\eqno(1.1)
$$
where $g_{ji}\in G$ is determined by
$$
h_jgh_i =g_{ji}h_k
$$
for some representative $h_k$.

\medskip

{\bf 1.2. Cohomology of $PSL(2,\bold{Z})$.} Let now
$G=PSL(2,\bold{Z})$ and $N$ a noncommutative $G$--module.
It is known that $PSL(2,\bold{Z})$ is the free product of
its two subgroups $\bold{Z}_2$ and $\bold{Z}_3$ generated respectively by
$$
\sigma = \left(\matrix 0& -1\\ 1& 0\endmatrix \right),\quad
\tau = \left(\matrix 0& -1\\ 1& -1\endmatrix \right)\, .
$$

$PSL(2,\bold{Z})$ acts transitively on $\bold{P}^1(\bold{Q}),$
the set of cusps of upper complex half--plane.
The stabilizer of $\infty$ is a cyclic subgroup $G_{\infty}$ generated by $\sigma\tau$. Hence the stabilizer $G_a$ of any cusp
$a\in \bold{P}^1(\bold{Q})$ is generated by $g^{-1}\sigma\tau g$
where $ga =\infty$.

\smallskip

Below we will give a concise description of the set $H^1(PSL(2,\bold{Z}),N)$
and its cuspidal subset $H^1(PSL(2,\bold{Z}),N)_{cusp}$
consisting by definition of those cocycle classes that become trivial
after restriction to any $G_a$.

\medskip

{\bf 1.2.1. Proposition.} {\it (i) Restriction to  $(\sigma ,\tau )$
of any cocycle in  $Z^1(PSL(2,\bold{Z}), N)$
belongs to the set
$$
\{\, (X,Y) \in N\times N\,|\, X\cdot\sigma X=1,\,
Y\cdot\tau Y\cdot\tau^2Y=1\,\}.
\eqno(1.2)
$$

\smallskip

(ii) Conversely, any element of the set (1.2) comes from
a unique 1--cocycle so that we can and will identify
these two sets. The cohomology 
relation between cocycles translates as
$$
(X,Y) \sim (n^{-1}X\sigma n, n^{-1}Y\tau n), \  n\in N.
\eqno(1.3)
$$

\smallskip

(iii) The cuspidal part of the cohomology consists of
classes of pairs of the form
$$
\{\,(X,Y)\,|\,\exists n\in N,\, X\cdot\sigma\tau Y=n^{-1}\cdot \sigma\tau n\,\}.
\eqno(1.4)
$$}

We will call (1.2) {\it abstract (noncommutative) Shimura--Eichler relations.} 

\smallskip

This result must be well known to experts, but I will
sketch a proof because I do not know a reference.

\smallskip

{\bf Proof.} Equations  (1.2) are a translation of the cocycle relations
applied to $\sigma^2=1$ and $\tau^3=1$.

\smallskip

Each nonidentical element $g$ of $PSL(2,\bold{Z})$
can be written uniquely as a product $g=\sigma^{a_1}\tau^{b_1} \dots
\sigma^{a_n}\tau^{b_n}$ with $n\ge 1$, $a_i=0$ or 1,
$b_i=0,1$ or 2, satisfying the condition that $a_i\ne 0$ for $i>1$
and $b_i\ne 0$ for $i<n$. Define {\it the length} $l(g)$ of such a
word as $\sum_i(a_i+b_i)$. Identity has length zero.

\smallskip

Each $g\ne 1$ ends with either $\sigma$, or $\tau$ that is,
can be represented as $h\sigma$ or $h\tau$ with $l(h)<l(g)$. 
All proofs proceed by induction on the length and use
the cocycle relations. Here are some details.
Denote by $Z$ the set (1.2) and by $\rho :\,Z^1(G,N)\to Z$
the restriction map.

\smallskip

(A) {\it $\rho$ is injective.} 

\smallskip

Assume that $X=u(\sigma ), Y=u(\tau )$ where $u$ is a cocycle.
Then we know $u$ on words of length $\le 1$. If
$g=h\sigma$ with $l(h)<l(g)\ge 2$, then
$u(g)=u(h)\cdot hX$ so that $u(g)$ is uniquely defined
by induction. The case $g=h\tau$ with $l(h)<l(g)$
is treated similarly. 

\smallskip

(B) {\it $\rho$ is surjective.}

\smallskip

Take arbitrary $(X,Y)\in Z.$ Construct a (well defined) map
$u:\,G\to N$ such that $u(1_G)=1_N, u(\sigma )=X$,
$u(\tau )=Y$ and $u(g)=u(h)\cdot hX$ (resp. $u(g)=u(h)\cdot hY)$)
if $g=h\sigma$ (resp. $g=h\tau$) and $l(h)<l(g).$

\smallskip

We have to check the cocycle relations (1.1): for arbitrary
$h,g\in G$, $u(hg)=u(h)\cdot hu(g).$ We make induction on $l(g)$.

\smallskip

{\it Start of induction: $l(g)=1$.} Then $g=\sigma$ or $\tau$,
and $l(hg)\ne l(h).$

\smallskip

If $l(hg)>l(h)$, then $hg$ ends with $g$, and the cocycle relation
holds by construction. 

\smallskip

If $l(hg)<l(h)$, then $h$ ends with $g$.
There are two subcases to consider: (a) $g=\sigma$,
$h=h^{\prime}\sigma$, $l(h^{\prime}) = l(h)-1$;
(b) $g=\tau$, $h=h^{\prime}\tau^2$, $l(h^{\prime})=l(h)-2.$
By construction, in the first case $u(h)=u(h^{\prime})\cdot h^{\prime}X$,
in the second case 
$$
u(h)=u(h^{\prime}\tau\cdot \tau )= u(h^{\prime}\tau)\cdot h^{\prime}\tau Y=
u(h^{\prime})\cdot h^{\prime}(Y\cdot \tau Y) =
u(h^{\prime})\cdot h^{\prime}((\tau^2Y)^{-1})=
u(h^{\prime})\cdot (hY)^{-1}.
$$
The relations we want to prove, namely
$u(h\sigma )=u(h)\cdot hX$ in the case (a), and
$u(h\tau )=u(h)\cdot hY$ in the case (b), easily follow.

\smallskip

{\it Inductive step: assuming cocycle relations
checked for all $g$ with $l(g)\le n-1$, check them for
longer words $g\sigma$ or $g\tau$ with $l(g)=n$.}

\medskip

In fact, applying the inductive assumption twice we get:
$$
u(hg\sigma )=u(hg)\cdot hgu(\sigma )= u(h)\cdot hu(g)\cdot hgu(\sigma )=
u(h)\cdot h[u(g)\cdot gu(\sigma )].
$$
On the other hand,
$$
u(h)\cdot hu(g\sigma )=u(h)\cdot h[u(g)\cdot gu(\sigma )].
$$
One treats $g\tau$ similarly. Thus $\rho$ is a bijection.

\smallskip

(C) The equivalence relations between cocycles clearly
restricts to (1.3) on $Z$. On the other hand,
if we start with a pair $(X,Y)$ and produce a cocycle $u$,
and then start with $(n^{-1}X\sigma n, n^{-1}Y\tau n)$
and produce another cocycle $v$, then by induction
one can check that $v(g)=n^{-1}u(g)\cdot gn$ for any $g$.
We leave this calculation to the reader.

\smallskip

(D) Finally, (1.4) means exactly that our cohomology class
becomes trivial on $G_{\infty}$, and triviality
on other  stabilizers of cusps follows from this.

\medskip

The description given above can be further cut down under some
additional conditions.

\medskip

{\bf 1.2.2. Proposition.} {\it Assume that the equation $n^2=m$ has a unique solution $n:=m^{1/2}$ in $N$ for any $n\in N.$ Then we have:

\smallskip
 
(i) Any cocycle is homologous to one
taking value 1 at $\sigma$.
Hence $H^1(PSL(2,\bold{Z}), N)$ can be identified with the following
quotient set:
$$
\{\, Y \in N\,|\, 
Y\cdot\tau Y\cdot\tau^2Y=1\,\}/ (Y\sim m^{-1}Y \tau m\
for\ some\ m\in N^{\sigma}).
\eqno(1.5)
$$

(ii) The cuspidal part of the cohomology consists of
classes of elements $Y$ of the form  $Y=n^{-1}\sigma\tau n$.}

\smallskip

{\bf Proof.}  Obviously $(gm)^{1/2}= g(m^{1/2})$
for all $g\in G,\, m\in N.$ 

\smallskip

Hence from
$X\cdot \sigma X=1$ it follows that $n^{-1}X \sigma n=1$
where $n=X^{1/2}$. Thus any cocycle is homologous to one with $X=1$.
On the subset of such cocycles, which we will identify with a part
of $N$ satisfying the second relation in (1.2),
the homology relation (1.3) becomes (1.5), and the cuspidal
relation (1.4) becomes (ii).

\smallskip

Whenever in $N$ equations $n^3=m$ are uniquely solvable {\it and
$N$ is commutative}, one can similarly show that
any cocycle (represented by) $(X,Y)$ is homologous to one with
$Y=1$: take $n=Y^{2/3} (\tau Y)^{1/3}$ in (1.3).
I was unable to check this in the noncommutative case.
However, in 2.6 we will see that  the iterated Shimura
cohomology class for $PSL(2,\bold{Z})$ can anyway be represented by a cocycle
with $Y=1$. Hence the following description parallel
to the one in 1.2.2 will be relevant;

\medskip

{\bf 1.2.3. Proposition.} {\it (i) The part of the set 
$H^1(PSL(2,\bold{Z}), N)$ represented by cocycles taking value 1
on $\tau$ can be identified with the following
quotient set:
$$
\{\, X \in N\,|\, 
X\cdot\sigma X =1\,\}/ (X\sim n^{-1}X \sigma n\
for\ some\ n\in N^{\tau}).
\eqno(1.6)
$$

(ii) The cuspidal part of the cohomology consists of
classes of elements $X$ of the form  $X=n^{-1}\cdot\sigma\tau n$.}

\bigskip

\centerline{\bf \S 2. Iterated Shimura integrals}

\medskip

{\bf 2.1. Forms of cusp modular type.} Let $\Gamma$ be 
a subgroup of finite index
of $SL(2,\bold{Z})$, $k\ge 2$ an integer, $S_k(\Gamma )$ the space
of cusp forms of weight $k$. Denote by $Sh_k(\Gamma )$ the space
of 1--forms on the complex upper half plane $H$ of the form
$f(z)\,P(z,1)\,dz$ where $f\in S_k(\Gamma )$, and $P=P(X,Y)$ runs over 
homogeneous polynomials
of degree $k-2$ in two variables.  Thus, the space $Sh_k(\Gamma )$
is spanned by  1--forms {\it of cusp modular type with integral
Mellin arguments in the critical strip} in the
terminology of [Ma3], Def. 2.1.1.

\medskip

{\bf 2.2. Action of $GL^+(2,\bold{R})$.} The group
of real matrices with a positive determinant $GL^+(2,\bold{R})$ acts on $H$
by fractional linear transformations $z\mapsto [g]z.$ Let
$j(g,z):=cz+d$ where $(c,d)$ is the lower row of $g$.
Then we have, for any function $f$ on $H$:
$$
g^*[f(z)\,P(z,1)\,dz] :=f([g]z) \,P([g]z,1)\,d([g]z)
$$
$$
=f([g]z)\,(j(g,z))^{-k}P(az+b,cz+d)\,\roman{det}\,g\,dz
\eqno(2.1)
$$  
where $(a,b)$ is the upper row of $g$. From the definition
it is clear that the diagonal matrices
act identically so that we have in fact an
action of $PGL^+(2,\bold{R})$. 

\smallskip

This can be rewritten in terms of the weight $k$ action
of $GL^+(2,\bold{R})$ upon functions on $H$. Actually,
in the literature one finds at least two different 
normalizations of such an action. They differ by a 
determinantal twist and therefore coincide
on $SL(2,\bold{R})$. For example, in [He1], [He2] one finds
$$
f\,|[g]_k(z):=f([g]z)\,j(g,z)^{-k}\,(\roman{det}\,g)^{k-1},
\eqno(2.2a)
$$
whereas in [Me2] and [Ma3] the action
$$
f\,|[g]_k^{\prime}(z):=f([g]z)\,j(g,z)^{-k}\,(\roman{det}\,g)^{k/2}
\eqno(2.2b)
$$
is used. 

\smallskip

Comparing this with (2.1), we get
$$
g^*[f(z)\,P(z,1)\,dz] =f\,|[g]_k(z)\,P(az+b,cz+d)\,(\roman{det}\,g)^{2-k}dz
$$
$$
=f\,|[g]_k^{\prime}(z)\,P(az+b,cz+d)\,(\roman{det}\,g)^{(2-k)/2}dz.
$$
Since $S_k(\Gamma )$ consists of holomorphic functions
which are $\Gamma$--invariant with respect to the (coinciding)
right actions (1.2a), (1.2b), the space $Sh_k(\Gamma )$ 
is $\Gamma$--stable and can be viewed as a tensor product
of the trivial representation on $S_k(\Gamma)$
and the $(k-2)$--th symmetric power of the basic
2--dimensional representation: for $g\in \Gamma$ we have
$$
g^*(f(z)\,P(z,1)\, dz)=f(z)\,P(az+b,cz+d)\,dz.
\eqno(2.3)
$$
\medskip

{\bf 2.3. Space $Sh_{\underline{k}}$ and formal series.}
In the following we choose and fix a group $\Gamma$ as above
and a finite family of pairwise distinct weights $\underline{k}=
(k_i).$ Put $Sh_{\underline{k}}:= \oplus_i Sh_{k_i}$. Denote
by $Sh_{\underline{k}}^*$ the dual space to $Sh_{\underline{k}}$,
together with the adjoint left action $g_*$ of $\Gamma$ on it,
so that $(g^*(\omega ),\nu )= (\omega ,g_*(\nu ))$ for all
$\omega \in Sh_{\underline{k}}$, $\nu \in Sh_{\underline{k}}^*$,
and $g\in \Gamma$.

\smallskip

We will consider the completed tensor algebra of 
$Sh_{\underline{k}}^*$ as a ring of formal series in a finite number
of associative non--commutative variables. Using the conventions
of [Ma3], we may and will choose a basis $(\omega_v)$ of $Sh_{\underline{k}}$
indexed by a finite set $V$, and the dual basis $(A_v)$
of $Sh_{\underline{k}}^*$. Then $\Gamma$ acts on the left
by linear transformations $g_*$ on $(A_v)$ inducing automorphisms
on the formal series ring $\bold{C}\langle\langle A_v\rangle\rangle.$ 
This ring has a continuous comultiplication defined by
$\Delta (A_v)=A_v\otimes 1+ 1\otimes A_v.$

\smallskip

{\it Group--like elements} $F$ of $\bold{C}\langle\langle A_v\rangle\rangle$
are characterized by the property $\Delta(F) = F\otimes F,\,
F\equiv 1\,\roman{mod}\,(A_v)$. As is well known, 
$F$ is group--like if and only if
$\roman{log}\,F$ belongs to the completed free Lie algebra freely generated
by $(A_v)$ inside $\bold{C}\langle\langle A_v\rangle\rangle$. 

\smallskip

We may extend the scalars $\bold{C}$ of
$\bold{C}\langle\langle A_v\rangle\rangle$ to functions or
1--forms on $H$. All scalars are assumed to commute with $(A_v)$.

\smallskip

In particular, the $\bold{C}\langle\langle A_v\rangle\rangle$--bimodule
$\Omega^1_H\langle\langle A_v\rangle\rangle$
contains a canonical element
$$
\Omega:=\sum_v A_v\omega_v \in Sh_{\underline{k}}^*\otimes Sh_{\underline{k}}
\eqno(2.4)
$$ 
which does not depend on the initial choice of basis $(\omega_v )$.

\medskip

{\bf 2.4. Iterated Shimura cocycles.} We will consider now
the iterated Shimura integrals
$$
J_a^z(\Omega):= 1+\sum_{n=1}^{\infty}
\int_{a}^z \Omega (z_1)
\int_{a}^{z_1} \Omega (z_2)\dots
\int_{a}^{z_{n-1}}\Omega (z_n)\,
\eqno(2.5)
$$
where $a,z$ are points of $\overline{H}:=H\,\cup\,\bold{P}^1(\bold{Q})$.
Such an integral is well defined and takes values in the group $\Pi$
of group--like elements of $\bold{C}\langle\langle A_v\rangle\rangle$.
For more details, see [Ma3]. 

\smallskip

The group $P\Gamma$ acts on $\Pi$ as was described in 1.3.

\smallskip

The following result is a slightly more precise
version of [Ma3], 2.6.1.

\medskip

{\bf 2.4.1. Theorem.} {\it (i) For any $a\in \overline{H}$, the map
$P\Gamma \to \Pi :\, \gamma\mapsto J^a_{\gamma a}(\Omega )$ is a noncommutative
1--cocycle in $Z^1(P\Gamma ,\Pi )$.

\smallskip

(ii) The cohomology class of this cocycle in $H^1(P\Gamma ,\Pi )$
does not depend on the choice of the reference point 
$a\in \overline{H}$.

\smallskip

(iii) This cohomology class belongs to the cuspidal subset
$H^1(P\Gamma ,\Pi )_{cusp}$ consisting of those
cohomology classes whose restrictions on the 
stabilizers of cusps in $\Gamma$
are trivial.}

\smallskip

The last statement which was not mentioned in [Ma3]
can be checked as follows. Let $\gamma$ belong to the stabilizer
$\Gamma_a$ of a point $a\in \bold{P}^1(\bold{Q}).$
Then if we take $a$ for the reference point,
the respective cocycle  is identically 1 on $\Gamma_a$
since $J_{\gamma a}^a=J_a^a=1$.

\medskip

{\bf 2.5. Reductions of the coefficient group.} We will call the
class $\zeta \in H^1(P\Gamma, \Pi )$, represented by $u_{\gamma}:=
J^a_{\gamma a}(\Omega )$
{\it the Shimura class}. The same name
will be applied to its various incarnations
obtained by changing $\Pi$ or $P\Gamma$ (and using Shapiro Lemma).

\smallskip

In this subsection, we will cut down $\Pi$
and exhibit representatives of this class
satisfying conditions stated in 1.2.2 and 1.2.3.

\medskip

{\bf 2.5.1. The continued fractions trick.} The following
result which we reproduce from [Ma3] drastically reduces
the size of a subgroup of $\Pi$ containing  a representative
of $\zeta$. 

\smallskip

Choose a set of representatives $C$ of left cosets
$P\Gamma \setminus PSL_2(\bold{Z})$. The iterated integrals
of the form $(J_{g(i\infty )}^{g(0)})^{\pm 1}$, $g\in C$, 
will be called {\it primitive} ones.
Notice that when $g\notin \Gamma$ the space spanned by
$(\omega_v)$ is not generally $g^*$--stable so that we cannot define
$g_*$.

\smallskip

{\bf 2.5.2. Proposition.} {\it Each $J^a_b(\Omega )$, $a,b\in\bold{P}^1(\bold{Q})$,
in particular components of any Shimura cocycle
with a cuspidal initial point $a\in \bold{P}^1(\bold{Q})$, can be expressed as
a noncommutative monomial in $\gamma_*(J^c_d(\Omega ))$ where
$\gamma$ runs over $\Gamma$ and $J^c_d(\Omega )$ runs over primitive integrals.}

\smallskip

{\bf Proof.} In fact, it suffices to express in this way $J_{i\infty}^a$
for $a>0$ (we omit $\Omega$ for brevity). Produce
a sequence of matrices $g_k$ from the consecutive convergents to $a$:
$$
a=\frac{p_n}{q_n},\ \frac{p_{n-1}}{q_{n-1}},\ \dots ,\ \frac{p_0}{q_0}=
\frac{p_0}{1},\ \frac{p_{-1}}{q_{-1}}:=\frac{1}{0},
$$
$$
g_k:=\left(\matrix p_{k} & (-1)^{k-1}p_{k-1}\\q_{k} & (-1)^{k-1}q_{k-1} \endmatrix
\right),\quad k=0, \dots , n.
$$
We have $g_k=g_k(a) \in SL_2(2,\bold{Z}).$
Put $g_k=\gamma_kc_k$ where $\gamma_k\in\Gamma$ and
$c_k\in C$ are two sequences of matrices depending on $a$.
Then
$$
J_{i\infty}^a=\prod_{k=n}^0 \gamma_{k*}(J^{c_k(i\infty)}_{c_k(0)} )
$$
which ends the proof.

\medskip

{\bf 2.6. The case $P\Gamma = PSL(2,\bold{Z}).$} In this case,
one can directly apply Proposition 1.2.1 providing noncommutative
Shimura-Eichler relations between iterated integrals. It can be also
applied to the smaller coefficient group $\Pi_0$ generated
by the $g_*(J_{i\infty}^0)$ as in 2.5.3.

\smallskip

In $\Pi$, any element has a unique square root.
Hence one can apply Proposition 1.2.2 as well. 
Another way to  produce a Shimura cocycle taking value 1
at $\sigma$, without using square roots, is to choose
$a=i$ for the initial point, because $\sigma i=i$.
Then the components of the respective Shimura cocycles
will be iterated integrals between points of complex multiplication
by $i$ rather than cusps, and the all--important $\tau$--component
is simply $Y=J_{\tau i}^i (\Omega ).$

\smallskip

A similar trick is applicable to $\tau$: its fixed point is
$\rho = e^{\pi i/3}$ and so we get $Y=1$,
$X=J_{\sigma\rho}^\rho (\Omega ).$

\smallskip

Notice finally that 
$$
J_{\sigma\rho}^\rho (\Omega ) = J_i^{\rho}(\Omega )\cdot \sigma
(J_i^{\rho}(\Omega ))^{-1}
$$
because $J_i^{\rho}=J_i^{\rho} J_{\sigma{\rho}}^i$
and $J_{\sigma{\rho}}^i=\sigma(J_{\rho}^i)=\sigma(J^{\rho}_i)^{-1}.$

\medskip

{\bf 2.7. Application of the Shapiro Lemma.} We can apply the Shapiro
Lemma for $P\Gamma \subset PSL(2,\bold{Z})$
in order to be able to use Proposition 1.2.1 for arbitrary
$\Gamma$.
The $P\Gamma$--module $\Pi$ gets replaced by the module $\Pi_{P\Gamma}$
of $P\Gamma$--covariant maps $PSL(2,\bold{Z})\to \Pi$.
Formula (1.1) shows that the Shimura class is still
represented by a cocycle whose components are iterated integrals
between two cusps. Square roots still exist and are unique
in $\Pi_{P\Gamma}$ so that Proposition 1.2.2 (i) is applicable as well.
However, the description of the cuspidal subset becomes somewhat
clumsier.  

\bigskip

\centerline{\bf \S 3. Linear term of $J_a^z(\Omega )$}

\smallskip

\centerline{\bf and classical modular symbols}

\medskip

The linear (in $(A_v)$) term of $J_a^z(\Omega )$ 
involves ordinary integrals of the form
$\int_a^z f(z)z^{s-1}dz$, $f\in S_k(\Gamma )$, $s\in \bold{C}$
(Mellin transform) or $s=1,\dots , k-1$ (Shimura integrals). 

\smallskip

In this section, I review some basic facts of the
classical theory of such integrals and explain how 
they extend (or otherwise) to the iterated setting,
following [Ma3].

\medskip

{\bf 3.1. Classical and iterated Mellin transforms.}
The classical Mellin transform of $f\in S_k(\Gamma )$
is
$$
\Lambda (f;s):=\int_{i\infty}^0f(z) z^{s-1}dz .
$$
Let $N>0$ and assume that $\Gamma$ is normalized by 
$$
g=g_N:=\left(\matrix 0& -1\\ N& 0\endmatrix \right) .
$$
Then $[N^{-1/2}g_N]_k$ defines an involution on $S_k(\Gamma)$
(see (2.2a)). Let $f$ be an eigenform with eigenvalue 
$\varepsilon_f=\pm 1$ with respect to this involution. Then
$$
\Lambda (f;s) = -\varepsilon_f e^{\pi is} N^{k/2-s} \Lambda (f; k-s)\,.
$$

\medskip

{\bf 3.1.1. The iterated extension.} The iterated Mellin transform of {\it a finite sequence of cusp forms}
$f_1,\dots ,f_k$ with respect to $\Gamma$ was defined in [Ma3]
as follows. Put $\omega_j(z):= f_j(z)\, z^{s_j-1} dz.$ Then
$$
M (f_1,\dots ,f_k;s_1,\dots ,s_k) :=
I_{i\infty}^0(\omega_1,\dots ,\omega_k)=
$$
$$
=\int_{i\infty}^{0} \omega_{1}(z_1)
\int_{i\infty}^{z_{1}} \omega_{2}(z_{2})\dots
\int_{i\infty}^{z_{n-1}}\omega_n(z_n)
$$
\smallskip
 
A neat functional equation however can be written not for these individual
integrals but for their generating series. More precisely,
let  $f_V=(f_v\,|\,v\in V)$ be a finite family
of cusp forms with respect to $\Gamma$,
$s_V=(s_v\,|\,v\in V)$ a finite family of complex numbers, $\omega_V=(\omega_v)$, where
$\omega_v(z):= f_v(z)\, z^{s_v-1} dz.$
The total Mellin transform 
of $f_V$ is
$$
TM(f_V;s_V) := J_{i\infty}^0(\omega_V) =
$$
$$
=1+\sum_{n=1}^{\infty}\sum_{(v_1,\dots ,v_n)\in V^n}
A_{v_1}\dots A_{v_n}\,M (f_{v_1},\dots ,f_{v_n};s_{v_1},\dots ,s_{v_n})
$$
Let $k_v$ be the weight of $f_v(z)$, and $k_V=(k_v)$.
Then we have
$$
TM(f_V;s_V)= g_{N*}(TM(f_V; k_V-s_V))^{-1}
$$
for an appropriate linear transformation $g_{N*}$
of formal variables $A_v$.

\medskip

{\bf 3.2. Dirichlet series.} It is well known that
$\Lambda (f;s)$ for general $s$ can be represented
by a product of a $\Gamma$--factor and a formal Dirichlet
series convergent in a right half plane of $s$:
$$
f(z) = \sum_{n=1}^{\infty} a_n e^{2\pi i nz}\quad \Longrightarrow\quad
\Lambda (f;s) = -\frac{\Gamma (s)}{(-2\pi i)^s} \sum_{n=1}^{\infty}
\frac{a_n}{n^s}
$$
In [Ma3], \S 3, it was shown that the iterated
Mellin transforms {\it at integral points of the product of critical strips} 
can be expressed as multiple Dirichlet series of a
special form. We omit the precise statements here.

\medskip

{\bf 3.3. The problem of iterated Hecke operators.} If $\Gamma$ is a congruence
subgroup, there is a well known classical correspondence
between the cusp forms which are eigenfunctions for the Hecke
algebra and their Mellin transforms admitting an Euler product.
Moreover:

\medskip

{\it (i) Shimura integrals of such a form span over $\overline{\bold{Q}}$
a linear space of dimension $\le 2$.}

\smallskip

This was proved in [Ma2] for $\Gamma = SL(2,\bold{Z})$,
and in [Sh3] for arbitrary (non necessarily congruence) $\Gamma$.

\smallskip

{\it (ii) In the case of a congruence subgroup $\Gamma$, Fourier 
coefficients of such forms are expressed by explicit
formulas involving summation of some simple linear functionals
over universal sets of matrices}.

\smallskip

This was also proved  in
[Ma2] for $\Gamma = SL(2,\bold{Z})$, and extended in several papers
to general congruence $\Gamma$. For an especially neat
version, see Merel's ``universal Fourier expansion'' in [Me2].

\medskip

The problem of extending these results to the iterated case
remains a major challenge. One obstacle is that
correspondences (in particular, Hecke correspondences)
do not act directly on the fundamental groupoid
(as opposed to the cohomology) and hence do not act
on the iterated integrals which
provide homomorphisms of this groupoid.

\smallskip

However, a part of the theory which is used in (ii),
that of the classical modular symbols,
allows a partial iterated extension. We will give below
a brief review of this theory.

\medskip

{\bf 3.4. Classical modular symbols.} The space of modular symbols
$MS_k(\Gamma )$, by definition, is essentially the space of {\it linear functionals}
on $S_k(\Gamma )$ spanned by the Shimura integrals
$$
f(z) \mapsto \int_{\alpha}^{\beta} f(z)z^{m-1}dz; \quad 1\le m \le k-1;
\quad \alpha, \beta \in \bold{P}^1(\bold{Q}).
$$
(but see more precise information below).
Three descriptions of $MS_k(\Gamma )$ are known:

\smallskip
 
(i) {\it Combinatorial (Shimura -- Eichler -- Manin)}: generators and relations.

\smallskip

(ii) {\it Geometric (Shokurov)}: $MS_k(\Gamma )$ can be identified
with a (part of) the middle homology of
the Kuga--Sato variety $M^{k)}$.

\smallskip

(iii) {\it Cohomological (Shimura)}: The dual space to $MS_k(\Gamma )$ can be identified with the cuspidal group cohomology $H^1(\Gamma , W_{k-2})_{cusp}$,
with coefficients in the $(k-2)$--th symmetric power of
the basic representation of $SL(2)$.

\smallskip

The noncommutative cohomology sets that we have described
in \S 2, are irerated extensions of this last description.

\smallskip

Here are some details.

\smallskip

{\bf 3.5. Combinatorial modular symbols.} In this description,
$MS_k(\Gamma )$ appears as an explicit subquotient of the space
$W_{k-2}\otimes \overline{C}$ where $W_{k-2}$ consists of polynomial forms
$P(X,Y)$ of degree $k-2$ of two variables, and $\overline{C}$
is the space of formal linear combinations of pairs of cusps
$\{\alpha,\beta\}\in \bold{P}^1(\bold{Q})$. Coeficients of these linear
combinations can be $\bold{Q}$, $\bold{R}$ or $\bold{C}$,
as in the theory of Hodge structure.

\smallskip

Each element of the form $P\otimes \{\alpha,\beta\}$ produces
a linear functional 
$$
f\mapsto \int_{\beta}^{\alpha} f(z)\,P(z,1)dz.
$$
This is extended to the total $W_{k-2}\otimes \overline{C}$ by linearity.

\smallskip

Denote by $C$ the quotient of $\overline{C}$ by the subspace generated
by sums $\{\alpha ,\beta\} + \{\beta ,\gamma\} + \{\gamma ,\alpha\}$.
Since $\int^\alpha_\beta +\int^\beta_\gamma + \int^\gamma_\alpha =0$,
our linear functional (Shimura integral) descends to $W_{k-2}\otimes C$. We will still
denote by $P\otimes \{\alpha,\beta\}$ the class of this element in $C$.

\smallskip

The group $GL^+(2,\bold{Q})$ acts from the left upon $W_{k-2}$
by $(gP)(X,Y):=P(bX-dY,-cX+aY)$ (notation as in (2.1)), and
upon $C$ by $g \{\alpha ,\beta\}:=\{g\alpha ,g\beta\}$.
Hence it acts on the tensor product.
A change of variable formula then shows that the Shimura
integral restricted to $S_k(\Gamma )$ vanishes on the subspace 
of $W_{k-2}\otimes C$ spanned  
by $P\otimes \{\alpha,\beta\} -gP\otimes \{g\alpha,g\beta\}$
for all $P\in W_{k-2}$, $g\in \Gamma$.

\smallskip

Denote by $MS_k(\Gamma )$ the quotient of  $W_{k-2}\otimes C$
by the latter subspace. 

\smallskip

The subspace of cuspidal modular symbols $MS_k(\Gamma )_{cusp}$ 
is defined by the following construction. Consider
the space $B$ freely spanned by $\bold{P}^1(\bold{Q})$.
Define the s pace $B_k(\Gamma )$ as the quotient of
$W_{k-2}\otimes B$ by the subspace generated
by $P\otimes \{\alpha\} - gP\otimes \{g\alpha\}$ for all
$g\in \Gamma$. There is a well defined boundary
map $MS_k(\Gamma )
\to B_k(\Gamma )$ induced by
$P\otimes \{\alpha ,\beta\} \mapsto P\otimes \{\alpha\} -
P\otimes \{\beta\}$. Its kernel is denoted 
$MS_k(\Gamma )_{cusp}$. 

\smallskip

By construction, any (real) modular
symbol in $MS_k(\Gamma )_{cusp}$ defines a $\bold{C}$--valued
functional $\int$ on $S_k(\Gamma )$ and in fact even
on $S_k(\Gamma ) \oplus \overline{S}_k(\Gamma )$.

\smallskip 
   
The first result of the theory is:

\medskip

{\bf 3.5.1. Theorem (Shimura).} {\it $\int$ is an isomorphism of $MS_k(\Gamma )_{cusp}$
with the dual space of $S_k(\Gamma ) \oplus \overline{S}_k(\Gamma )$.}  

\medskip

{\bf 3.5.2. Remark.} Probably, the most interesting recent result 
involving combinatorial modular symbols is Herremans'
combinatorial reformulation of Serre's conjecture in [He1], [He2].

\medskip

{\bf 3.6. Geometric modular symbols.} Let $\Gamma^{(k)}$ be
the semidirect product
$\Gamma\ltimes (\bold{Z}^{k-2}\times \bold{Z}^{k-2})$
acting upon $H\times \bold{C}^{k-2}$ via
$$
(\gamma ;\,n,m)\,(z,\zeta ):= ([\gamma ]z;\,j(\gamma ,z)^{-1}
(\zeta +zn +m))
$$
where $n=(n_1,\dots ,n_{k-2})$, 
$m=(m_1,\dots ,m_{k-2})$, $\zeta =(\zeta_1,\dots ,\zeta_{k-2})$,
and $nz=(n_1z,\dots ,n_{k-2}z)$.

\smallskip

If $f(z)$ is a $\Gamma$--invariant cusp form of weight $k$, then
$$
f(z)dz\wedge d\zeta_1\wedge \dots \wedge d\zeta_{k-2}
$$
is a $\Gamma^{(k)}$--invariant holomorphic volume form
on $H\times \bold{C}^{k-2}$. Hence one can push it down
to a Zariski open smooth subset of the quotient
$\Gamma^{(k)} \setminus (H\times \bold{C}^{k-2})$.
An appropriate smooth compactification $M^{(k)}$ of this subset
is called {\it a Kuga--Sato variety,} cf. [Sh1]--[Sh3].

\smallskip

Denote by $\omega_f$ the image of this form on $M^{(k)}$.
Notice that it depends only on $f$, not on any Mellin argument.
The latter can be accomodated in the structure of (relative) cycles
in $M^{(k)}$, so that integrating $\omega_f$ over such cycles
we obtain the respective Shimura integrals.

\smallskip

Concretely, let $\alpha ,\beta \in \bold{P}^1(\bold{Q})$ 
be two cusps in $\overline{H}$ and let $p$ be a geodesic joining
$\alpha$ to $\beta$. Fix $(n_i)$ and $(m_i)$ as above.
Construct a cubic singular cell $p\times (0,1)^{k-2}\to H\times\bold{C}^{k-2}$:
$(z, (t_i))\mapsto (z, (t_i(zn_i+m_i)))$. Take the
$S_{k-2}$--symmetrization of this cell and
push down the result to the Kuga--Sato variety.
We will get a relative (modulo fibers of $M^{(k)}$ over cusps)
cycle whose homology class is Shokurov's
higher modular symbol $\{\,\alpha, \beta ;\,n,m\}_{\Gamma}.$
One easily sees that
$$
\int_{\alpha}^{\beta} f(z) \prod_{i=1}^{k-2}(n_iz+m_i)\,dz=
\int_{\{\,\alpha, \beta ;\,n,m\}_{\Gamma}} \omega_f\,.
$$
The singular cube $(0,1)^{k-2}$ may also be replaced by
an evident singular simplex.

\medskip

{\bf 3.6.1. Theorem (Shokurov).} {\it (i) The map $f\mapsto \omega_f$
is an isomorphism $S_k(\Gamma )\to H^0(M^{(k)}, \Omega_{M^{(k)}}^{k-1}).$

\smallskip

(ii) The homology subspace spanned by Shokurov modular symbols
with vanishing boundary
is canonically isomorphic to the space of cuspidal combinatorial modular symbols.}

\medskip

{\bf 3.6.2. Remark.} I suggested in [Ma3] that it would be desirable 
to replace in this description Kuga--Sato varieties by moduli spaces of 
curves of genus 1 with marked points and a level structure.
For $\Gamma = SL(2,\bold{Z})$, this was essentially accomplished in a recent
paper [CF] by C.~Consani and C.~Faber. Namely, they proved
that the Chow motive associated with $S_k(SL(2,\bold{Z}))$
(with coefficients in $\bold{Q}$) is cut off
by the alternating projector from the motive of $\overline{M}_{1,k-2}$.
Recall that the symmetric group $S_{k-2}$ renumbering
marked points naturally acts on  $\overline{M}_{1,k-2}$.

\medskip

{\bf 3.7. Cohomological modular symbols.} In this description, the space
dual to $MS_k(\Gamma )$ is identified with the group
cohomology $H^1(\Gamma , W_{k-2}).$

\smallskip

A bridge between the geometric and the cohomological 
descriptions is furnished by the identification of 
$H^1(\Gamma , W_{k-2})_{cusp}$
with the cohomology of a local system on $M_{1,1}$,
namely $H^1_!(M_{1,1}, \roman{Sym}^{k-2} R^1\pi_*\bold{Q})$.

\smallskip

Our iterated version explained in \S 2 was an attempt to
extend this version of modular symbols.

\medskip

{\bf 3.7.1. Remark.} $SL(2)$--modular symbols
(and their generalization to groups of higher rank)
made their appearance also in
the context of (relations between) multiple polylogarithms:
see A. Goncharov's papers [Go3], [Go4]. It is not clear (at least to me)
how to connect this description with the former ones.

\bigskip

\centerline{\bf Bibliography}

\medskip

[Ch] K.--T. Chen. {\it Iterated path integrals.}
Bull. AMS, 83 (1977), 831--879.

\smallskip

[CoF] C.~Consani, C.~Faber. {\it On the cusp form motives in
genus 1 and level 1.} e--Print math.AG/0504418 .

\smallskip

[DeGo]  P.~Deligne A.~Goncharov. {\it Groupes fondamentaux
motiviques de Tate mixte.} e--Print math.NT/0302267 .

\smallskip

[Dr1] V.~G.~Drinfeld. {\it Two theorems on modular curves.} Func. An.
and its Applications, 7:2 (1973), 155--156.

\smallskip

[Dr2] V.~G.~Drinfeld. {\it On quasi--triangular quasi--Hopf
algebras and some groups closely associated
with $Gal\,(\overline{\bold{Q}}/\bold{Q})$.}
Algebra and Analysis 2:4 (1990); Leningrad Math. J.
2:4 (1991), 829--860.

\smallskip

[Go1] A.~Goncharov. {\it Polylogarithms in arithmetic and geometry.}
Proc. of the ICM 1994 in Z\"urich,  vol. 1, Birkh\"auser, 1995, 374--387.

\smallskip

[Go2] A.~Goncharov. {\it Multiple $\zeta$-values, Galois
groups and geometry of
modular varieties}. Proc. of the Third European Congress of
Mathematicians. Progr. in Math.,
vol. 201, p. 361--392. Birkh\"auser Verlag, 2001.
e--Print math.AG/0005069.

\smallskip

[Go3] A.~Goncharov. {\it The double logarithm and Manin's complex for modular curves.} Math. Res. Letters, 4 (1997), 617--636.

\smallskip

[Go4] A.~Goncharov. {\it Multiple polylogarithms,
cyclotomy, and modular complexes.}  Math. Res. Letters, 5 (1998), 497--516.

\smallskip

[Go5] A.~Goncharov. {\it Multiple polylogarithms and mixed Tate motives.}
e--Print math.AG/0103059

\smallskip

[GoMa] A.~Goncharov, Yu.~Manin. {\it Multiple zeta-motives and moduli spaces $\overline{M}_{0,n}$.} Compos. Math.
140:1 (2004), 1--14.  e--Print math.AG/0204102

\smallskip

[Ha] R.~Hain. {\it Iterated integrals and algebraic cycles:
examples and prospects.} e--Print math.AG/0109204

\smallskip

[He1] A.~Herremans. {\it A combinatorial interpretation of 
Serre's conjecture on modular Galois representations.}
Thesis, Katholieke Universiteit Leuven, 2001.

\smallskip

[He2] A.~Herremans. {\it A combinatorial interpretation of Serre's conjecture
on modular Galois representations.} 
Ann. Inst. Fourier, 53:5 (2003), 1287--1321.

\smallskip

[Ma1] Yu.~Manin. {\it Parabolic points and zeta-functions of modular curves.}
 Russian: Izv. AN SSSR, ser. mat. 36:1 (1972), 19--66. English:
Math. USSR Izvestija, publ. by AMS, vol. 6, No. 1 (1972), 19--64,
and Selected papers, World Scientific, 1996, 202--247.

\smallskip

[Ma2] Yu.~Manin. {\it Periods of parabolic forms and $p$--adic Hecke series.}
Russian: Mat. Sbornik, 92:3 (1973), 378--401. English:
Math. USSR Sbornik, 21:3 (1973), 371--393
and Selected papers, World Scientific, 1996, 268--290.

\smallskip

[Ma3] Yu.~Manin. {\it Iterated integrals of modular forms and noncommutative modular symbols.} 37 pp. e--Print math.NT/0502576

\smallskip

[Me1] L.~Merel. {\it Quelques aspects arithm\'etiques et g\'eom\'etriques
de la th\'eorie des symboles modulaires.} Th\'ese de doctorat,
Universit\'e Paris VI, 1993.

\smallskip

[Me2] L.~Merel. {\it Universal Fourier expansions of modular forms.}
Springer Lecture Notes in Math., vol. 1585 (1994), 59--94

\smallskip
 
[PlRap] V.~Platonov, A.~Rapinchuk. {\it Algebraic groups and number
theory.} Pure and Applied Math., 139. Academic Press,
Boston, 1994

\smallskip

[Sh1] V.~Shokurov. {\it Modular symbols of arbitrary weight.} Func. An.
and its Applications, 10:1 (1976), 85--86.

\smallskip

[Sh2] V.~Shokurov. {\it The study of the homology of Kuga
varieties.} Math. USSR Izvestiya, 16:2 (1981), 399--418.

\smallskip

[Sh3] V.~Shokurov. {\it Shimura integrals of cusp forms.}
Math. USSR Izvestiya, 16:3 (1981), 603--646.

\smallskip

[Za1] D.~Zagier. {\it Hecke operators and periods of modular forms.}
In: Israel Math. Conf. Proc., vol. 3 (1990), I.~I.~Piatetski--Shapiro
Festschrift, Part II, 321--336.

\smallskip

[Za2] D.~Zagier. {\it Values of zeta functions and their
applications.} First European Congress of Mathematics,
vol. II, Birkh\"auser Verlag, Basel, 1994, 497--512.

\smallskip

[Za3] D.~Zagier. {\it Periods of modular forms and Jacobi
theta functions.} Inv. Math., 104 (1991), 449--465.

\enddocument